\documentclass[10pt]{article}
\usepackage{amssymb,amsmath,graphicx,acronym,amsfonts}
\usepackage{color}

\newtheorem{corollary}{Corollary}[section]
\newtheorem{theorem}{Theorem}[section]
\newtheorem{lemma}{Lemma}[section]
\newtheorem{definition}{Definition}[section]
\newtheorem{proposition}{Proposition}[section]
\newtheorem{example}{Example}[section]
\newtheorem{assum}{Assumption}[section]
\newtheorem{algo}{Algorithm}[section]
\newtheorem{Remark}{Remark}[section]

\def\bc{\begin{corl}}
\def\bc{\end{corl}}
\def\ba{\begin{algo}}
\def\ea{\end{algo}}
\def\br{\begin{Remark}}
\def\er{\end{Remark}}
\def\bs{\begin{assum}}
\def\es{\end{assum}}
\def\bt{\begin{theorem}}
\def\et{\end{theorem}\vskip 3pt}
\def\bl{\begin{lemma}}
\def\el{\end{lemma}}
\def\ep{\end{proposition}}
\def\bp{\begin{proposition}}
\def\qed{\hfill{$\Box$}\vskip 5pt}
\def\be{\begin{example}}
\def\ee{\end{example}}
\def\bd{\begin{definition}}
\def\ed{\end{definition}}
\def\bc{\begin{corollary}}
\def\ec{\end{corollary}}
\def\proof{\noindent\it Proof. \hspace{1mm}\rm}
\input epsf
\begin{document}
\title{\bf Further Results on Cauchy Tensors and Hankel Tensors}
\author{Haibin Chen\thanks{Department of Applied Mathematics, The Hong Kong Polytechnic University, Hung Hom,
Kowloon, Hong Kong. Email: chenhaibin508@163.com. This author's work was supported by the Natural Science Foundation of China (11171180).} \quad
Guoyin Li \thanks{Department of Applied Mathematics, University of New South Wales, Sydney 2052, Australia. E-mail: g.li@unsw.edu.au (G.
Li). This author's work was partially supported by Australian Research Council.} \quad
Liqun Qi
\thanks{Department of Applied Mathematics, The Hong Kong Polytechnic University, Hung Hom,
Kowloon, Hong Kong. Email: maqilq@polyu.edu.hk. This author's work
was supported by the Hong Kong Research Grant Council (Grant No.
PolyU 502111, 501212, 501913 and 15302114).} }
\maketitle
\vspace{-0.3cm}
\begin{abstract} {  In this article, we present various new results on Cauchy tensors and Hankel tensors.
{  We first introduce the concept of generalized Cauchy tensors which extends Cauchy tensors in the current literature, and provide several conditions characterizing positive semi-definiteness of generalized Cauchy tensors with nonzero entries.}
As a consequence, we show that Cauchy tensors are positive semi-definite if and only if they
are SOS (Sum-of-squares) tensors.} Furthermore, we prove that all positive semi-definite Cauchy tensors are
completely positive tensors, which means every positive semi-definite Cauchy tensor can be decomposed {  as} the sum of
nonnegative rank-1 tensors. We also establish that all the H-eigenvalues of nonnegative Cauchy tensors are nonnegative.
Secondly, we present new mathematical properties of Hankel tensors.
{  We prove that an even order Hankel tensor is Vandermonde positive semi-definite if and only if its associated plane tensor is positive semi-definite.  We also show that,
if the Vandermonde rank of a Hankel tensor $\mathcal{A}$ is less than the dimension of the underlying space, then positive semi-definiteness
of $\mathcal{A}$ is equivalent to the fact that $\mathcal{A}$ is a complete Hankel tensor, and so, is further equivalent to the SOS property of $\mathcal{A}$.  Lastly, we introduce a new structured tensor called Cauchy-Hankel tensors,
which is a special case of Cauchy tensors and Hankel tensors simultaneously.}
Sufficient and necessary conditions are established for an even order Cauchy-Hankel tensor to be positive definite.  Final remarks are listed at the end of the paper.

\end{abstract}

\newpage
\section{Introduction}

Let $\mathbb{R}^n$ be the $n$ dimensional real Euclidean space.
Denote the sets of all natural numbers by $\mathbb{N}$. Suppose $m, n\in \mathbb{N},~m, n\geq2$
and denote $[n]=\{1,2,\cdots,n\}$.
It should be noted in advance, we always consider real order $m$
dimension $n$ tensors in this paper.

Tensors (or sometimes called hypermatrices) are the multi-array extensions of matrices.
It was recently demonstrated in \cite{Hillar13} that most of the problems associated with tensors are, in general, NP-hard.
So, it motivates researchers to study tensors with special structure i.e. structured tensors.
In the last two or three years, a lot of
research papers on structured tensors appeared
\cite{chen14,Chen14,Ding13,Ding14,HH14,LQX14,Qi14,qi14,QXX,
Song14,song14,YY14,Zhang12}. These include M tensors, circulant
tensors, completely positive tensors, Hankel tensors, Hilbert
tensors, P tensors, B tensors and Cauchy tensors.
Many interesting properties and meaningful results
of structured tensors have been discovered. For instance, spectral properties of structure tensors, positive
definiteness and semi-definiteness of structured tensors were established.
Furthermore, some practical applications of structured tensors were studied such as application in stochastic
process and data fitting \cite{Chen14, Ding14}. Very recently, authors of \cite{LQX14} studied
SOS-Hankel tensors and applied them to the positive semi-definite tensor completion problem.

{  Among the various structured tensors we mentioned above, there are two particular interesting classes: Cauchy tensors and Hankel tensors.
The symmetric Cauchy tensors were defined and analyzed in \cite{chen14}.
In the following discussion, {  we simply refer it as Cauchy tensors instead of symmetric Cauchy tensors.} One of the nice properties of a Cauchy
tensor is that its positive semi-definiteness (or positive definiteness) can be easily verified by the sign of the associated generating vectors.}
In fact, it was proved in \cite{chen14} that an even order Cauchy tensor is
positive semi-definite if and only if each of entries of its generating vector is positive, and an even
order Cauchy tensor is positive definite if and only if each entries of its generating vector
is positive and mutually distinct.

{  Hankel tensors arise from signal processing and data fitting \cite{BB, Ding14, Papy05}. As far as we know, the definition of Hankel tensor was first
introduced in \cite{Papy05}. Recently,  some easily verifiable structured tensors related to Hankel tensors were also introduced in \cite{Qi14}. These structured tensors include
strong Hankel tensors, complete Hankel tensors and the
associated plane tensors that correspond to underlying Hankel tensors.}  It was proved that if a Hankel tensor
is copositive or an even order Hankel tensor is positive semi-definite, then the associated plane tensor
is copositive or positive semi-definite respectively \cite{Qi14}. Furthermore, results on positive semi-definiteness
of even order strong and complete Hankel tensors were given. However, the relationship between strong Hankel tensors and
complete Hankel tensors was not provided in \cite{Qi14}. Later, in \cite{LQX14}, it was shown that complete Hankel tensors are strong Hankel tensors; while the converse
is, in general, not true.

In this paper, we will provide further new  results for
Cauchy tensors and Hankel tensors which complements the existing literature.
The remainder of this paper is organized as follows.
In Section 2, we will first recall some basic notions of tensors are given such as H-eigenvalues, Z-eigenvalues and
positive semi-definite tensors. We will also introduce the notion of Vandermonde positive semi-definite tensors, which is a special class of positive semi-definite tensors.

In Section 3, we will first introduce the generalized Cauchy tensors which is an extension of the Cauchy tensors in the literature.
Then, we provide complete characterization for positive semi-definite generalized Cauchy tensors with nonzero entries.
We will also present sufficient and necessary conditions guaranteeing an even order generalized Cauchy tensor with
nonzero entries to be completely positive. After that, similar conclusions on even order Cauchy tensors are
established. It is proven that all positive semi-definite Cauchy tensors are SOS (sum-of-square) tensors and an even order
Cauchy tensor is positive semi-definite if and only if it is completely positive tensor, which means every positive semi-definite
Cauchy tensor can be written as a sum of nonnegative rank-1 tensors.
Furthermore, we prove that the Hadamard product of two positive semi-definite Cauchy tensors are still positive semi-definite tensors. And the nonnegativity
for H-eigenvalues of nonnegative Cauchy tensors are testified.

In Section 4, we provide further new properties of Hankel tensors. We prove that the associated plane tensor of an even order Hankel tensor
is positive semi-definite if and only if the Hankel tensor is Vandermonde positive semi-definite. Using this conclusion,
we give an example to show that, for higher dimensional Hankel tensors, the associated plane tensor is positive semi-definite but the Hankel tensor failed.
Combining this with the methods in \cite{QWW}, one may easily get that the positive semi-definiteness of the associated plane tensor is a checkable necessary condition for the positive semi-definiteness of Hankel tensors. {
Sufficient conditions for a complete Hankel tensor to
be positive semi-definite Hankel tensor are also provided. }

In Section 5, we introduce Cauchy-Hankel tensors, which are natural extensions
of Cauchy-Hankel matrices. {  The class of Cauchy-Hankel tensors is a subset
of Cauchy tensors \cite{chen14} and Hankel tensors \cite{Ding14, LQX14, Qi14} simultaneously. We provide a checkable sufficient and necessary
condition for an even order Cauchy-Hankel tensor to be positive definite.
We also show that an even order Cauchy-Hankel tensor is positive
semi-definite if and only if the associated homogeneous polynomial is strict monotonically increasing
on the nonnegative orthant $\mathbb{R}_+^n$.}

Some final remarks are provided in Section 6.

{  Before we end the introduction section, let us make some comments on the symbols
that will be used throughout this paper. }
Vectors are denoted by italic lowercase letters i.e. $x,~ y,\cdots$,
and matrices are denoted by capital letters $A,~B,\cdots$.
Suppose $e\in \mathbb{R}^n$ be
all one vectors and let $e_i$ denotes the $i$-th unite coordinate vector in $\mathbb{R}^n$.
We use bold letters ${\bf 0} \in \mathbb{R}^n$ to denote zero vector.
Tensors are written as calligraphic capitals such as
$\mathcal{A}, \mathcal{T}, \cdots.$  Let $\mathcal{I}$
denote the real identity tensor. For $x=(x_1,x_2,\cdots,x_n)^T, y=(y_1,y_2,\cdots,y_n)^T\in \mathbb{R}^n$,
then $x\geq y$ ($x\leq y$) means $x_i\geq y_i$ ($x_i\leq y_i$) for all $i\in [n]$. $x^{[m]}$ is defined by
$(x_1^m,x_2^m,\cdots,x_n^m)^T$.

\setcounter{equation}{0}
\section{Preliminaries}
A real tensor with order $m$ and dimension $n$ is defined
by $\mathcal{A}=(a_{i_1i_2\cdots i_m})$, $i_j\in [n]$, $j\in [m]$. If the entries $a_{i_1i_2\cdots i_m}$ are
invariant under any permutation of the subscripts, then tensor $\mathcal{A}$ is called symmetric tensor. Let
$x=(x_1,x_2,\cdots,x_n)^T\in \mathbb{R}^n$. The two forms below will be used in the following analysis frequently:
$$\mathcal{A}x^{m-1}=\left(\sum_{i_2,i_3,\cdots,i_m=1}^na_{ii_2\cdots i_m}x_{i_2}\cdots x_{i_m}\right)_{i=1}^n;$$
$$\mathcal{A}x^m=\sum_{i_1,i_2,\cdots,i_m=1}^na_{i_1i_2\cdots i_m}x_{i_1}x_{i_2}\cdots x_{i_m}.$$

Denote $\mathbb{R}^n_+=\{x\in \mathbb{R}^n ~|~x \geq {\bf 0}\}$. If $\mathcal{A}x^m\geq 0$ for all $x\in \mathbb{R}^n_+$,
then $\mathcal{A}$ is called copositive.
An even order m dimension n tensor $\mathcal{A}$ is called positive semi-definite if for any vector $x\in \mathbb{R}^n$,
it satisfies $\mathcal{A}x^m\geq 0$. Tensor $\mathcal{A}$ is called positive definite if $\mathcal{A}x^m > 0$ for all
nonzero vectors $x\in \mathbb{R}^n$. {  From the definition, it is easy to see that, for a positive semi-definite tensor, its order $m$ must be an even number.
Therefore, in the following analysis, we always assume the order of the tensor is even when we consider a positive semi-definite tensor.}

We call $u\in \mathbb{R}^n$ a Vandermonde vector if  $u=(1, \mu, \mu^2,\cdots, \mu^{n-1})^T\in \mathbb{R}^n$ for some $\mu\in \mathbb{R}$.
If $\mathcal{A}u^m\geq 0$ for all Vandermonde vectors $u\in \mathbb{R}^n$,
then we say that tensor $\mathcal{A}$ is {\bf Vandermonde positive semi-definite}. It's obvious that
positive semi-definite tensors are always Vandermonde positive semi-definite, but not vice versa.

{  Next, we recall the definitions of eigenvalues of tensors.

\begin{definition}
Let $\mathcal {A}$ be a symmetric tensor with order $m$ and dimension $n$.
We say $\lambda \in \mathbb{R}$ is a $Z$-eigenvalue of $\mathcal {A}$
and $x \in \mathbb{R}^n \backslash \{\bf 0\}$ is an $Z$-eigenvector corresponding to $\lambda$  if $(x,\lambda)$ satisfies
\[
\left\{\begin{array}{c}
\mathcal{A}x^{m-1}=\lambda x,  \\
x^Tx=1.
       \end{array}\right.
\]
Moreover, we say $\lambda \in \mathbb{R}$ is an $H$-eigenvalue of $\mathcal {A}$
and $x \in \mathbb{R}^n \backslash \{ \bf 0\}$ is an $H$-eigenvector corresponding to $\lambda$  if $(x,\lambda)$ satisfies
\[
\mathcal{A}x^{m-1}=\lambda x^{[m-1]},
\]
where $x^{[m-1]}=(x_1^{m-1},\ldots,x_n^{m-1})^T \in \mathbb{R}^n$.
\end{definition}

The definitions of $Z$-eigenvalues and $H$-eigenvalues  were introduced by Qi in \cite{Qi05}.
Independently, Lim \cite{Lim05} also gave the definitions via a
variational approach and established an interesting Perron-Frobenius theorem for tensors with nonnegative entries. From \cite{Qi05} and \cite{CPZ09},
both $Z$-eigenvalues and $H$-eigenvalues for an even order symmetric tensor always exist. Moreover, from the definitions, we can see that
 finding an $H$-eigenvalue of a symmetric tensor is equivalent to solving a  homogeneous polynomial equation while calculating a $Z$-eigenvalue is equivalent to solving nonhomogeneous polynomial equations. In general, the behaviors of $Z$-eigenvalues and $H$-eigenvalues can be quite different. For example, a diagonal symmetric tensor A has exactly $n$ many $H$-eigenvalues and may have more than $n$ Z-eigenvalues (for more details see \cite{Qi05}).
 Recently, a lot of researchers have devoted themselves to the study of
 eigenvalue problems of symmetric tensors and have found important applications in diverse areas including spectral hypergraph theory \cite{Cooper,LQY2}, dynamical control \cite{Qi_Ni},
 medical image science \cite{LQY1,QYW} and signal processing \cite{La_book}.}

%
%
%
%
%
\setcounter{equation}{0}
\section{SOS properties and {  complete positivity} of even order Cauchy tensors}

Symmetric Cauchy tensors was first studied in \cite{chen14}. Some checkable sufficient and necessary conditions for
an even order symmetric Cauchy tensor to be positive semi-definite or positive definite tensors {  were}  provided in \cite{chen14}, which extends the matrix cases established in
\cite{Fied10}.
\begin{definition}\label{def31}$^{\cite{chen14}}$ Let $c=(c_1,c_2,\cdots,c_n)^T\in \mathbb{R}^n$. Let a real tensor $\mathcal{C}=(c_{i_1i_2\cdots i_m})$ be defined by
\begin{equation}\label{e31}
c_{i_1i_2\cdots i_m}=\frac{1}{c_{i_1}+c_{i_2}+\cdots+c_{i_m}},\quad j\in [m],~i_j \in [n].\end{equation}
Then, we say that $\mathcal{C}$ is a symmetric {\bf Cauchy tensor} with order $m$ and dimension $n$. The corresponding vector $c\in \mathbb{R}^n$ is called the generating vector of $\mathcal{C}$.
\end{definition}

Now, given two vectors $c=(c_1,c_2,\cdots,c_n)^T,~d=(d_1,d_2,\ldots,d_n)^T\in \mathbb{R}^n$. Consider the
{\bf generalized Cauchy tensor} $\mathcal{C}=(c_{i_1i_2\cdots i_m})$ with order $m$ dimension $n$, where
$$
c_{i_1 i_2\cdots i_m}=\frac{d_{i_1} d_{i_2} \cdots d_{i_m}}{c_{i_1}+c_{i_2}+\cdots+c_{i_m}}, \ i_j\in [n], \ j\in [m].
$$
For the sake of simplicity, we call vectors $c,~d$ generating vectors of
 generalized Cauchy tensor $\mathcal{C}$. In the special case when $d_i=1$, $i\in [n]$, a generalized
Cauchy tensor reduces to a Cauchy tensor defined in Definition \ref{def31}. In the case when $m=2$, a generalized Cauchy tensor collapses to a symmetric
generalized Cauchy matrix \cite{Polya}. We also note that every rank-one tensor with the form $u^ m$ for some $u \in \mathbb{R}^n$ is, in particular, a
generalized Cauchy tensor.

Define Cauchy tensor $\overline{\mathcal{C}}=(\overline{c}_{i_1, i_2,\cdots i_m})$ where
$$
\overline{c}_{i_1, i_2,\cdots i_m}=\frac{1}{c_{i_1}+c_{i_2}+\cdots+c_{i_m}}, \ i_j=1,\cdots,n, \ j=1,\cdots,m.
$$
It is easy to see for any $x \in \mathbb{R}^n$, we have
$$\mathcal{C} x^m \equiv \overline{\mathcal{C}} y^m,$$
where $y \in \mathbb{R}^n$ with $y_i = d_ix_i$ for $i = 1,\cdots,n$. By
Theorems 2.1 and 2.3 of \cite{chen14}, one may easily conclude that the
generalized Cauchy tensor $\mathcal{C}$ is positive semi-definite if
and only if $d_i=0, c_i\neq0$ or $d_i\neq0, c_i>0$, $i\in [n]$ and $\mathcal{C}$ is positive definite if and
only if $c_1, c_2, \cdots, c_n$ are positive {  real number} and mutually distinct,
and $d_i \neq 0$, $i=1,\cdots,n$.

In this section, we mainly characterize SOS (sum-of-squares) properties and completely positiveness of even order generalized Cauchy tensors
with nonzero entries. Then, similar results on even order Cauchy tensors are established.
Before giving the main results, we briefly {  recall the definitions of SOS tensors and completely positive tensors}.

SOS tensors are first defined in \cite{Hu14}. The definition of SOS tensors relies on the celebrated concept of SOS polynomials, which is a fundamental
concept in polynomial optimization theory \cite{HLQS,Hu14,JB01, M09, shor98}. Assume $\mathcal{A}$ is an order $m$ dimension
$n$ symmetric tensor. {  Let $m=2k$ be an even number. If $$f(x)=\mathcal{A}x^m,~x\in \mathbb{R}^n$$
can be decomposed to the sum of squares of polynomials of degree $k$, then $f$ is called a {\bf sum-of-squares (SOS) polynomial}, and the corresponding symmetric tensor $\mathcal{A}$ is called an {\bf SOS tensor} \cite{Hu14}.    From the definition, any SOS tensor is positive semi-definite. On
the other hand, the converse is not true, in general \cite{HLQS,Hu14}. The importance of studying SOS tensors is that the problem for determining an even order symmetric tensor is an SOS tensor or not is equivalent to solving a semi-infinite linear programming problem, which can be done in polynomial time;
while determining the positive semi-definiteness of a symmetric tensor is, in general, NP-hard. Interestingly, it was recently shown in \cite{Hu14}  that for a so-called Z-tensor $\mathcal{A}$ where the off-diagonal elements are all non-positive,
$\mathcal{A}$ is positive semi-definite if and only if it is a SOS tensor.}

{  Tensor $\mathcal{A}$ is called a {\bf completely decomposable tensor} if there are
vectors $x_j\in \mathbb{R}^n$,~$j\in [r]$ such that $\mathcal{A}$ can be written as sums of rank-one tensors generated by the vector $x_j$, that is,
$$\mathcal{A}=\sum\limits_{j\in [r]}x_j^m.$$
If $x_j\in \mathbb{R}^n_+$ for all $j\in [r]$, then $\mathcal{A}$ is called a {\bf completely positive tensor} \cite{QXX}. It was shown that a strongly
symmetric, hierarchically dominated nonnegative tensor is a completely positive tensor \cite{QXX}. It can be directly verified that all even order completely
positive tensors are SOS tensors, and so, are also positive semi-definite tensors. We note that verifying a tensor $\mathcal{A}$ is a completely decomposable or not, and
finding its explicit rank one decomposition are  highly nontrivial. This topic has attracted a lot of researchers and many
important work has been established along this direction. For detailed discussions, see \cite{Lim,Kolda,QXX} and the reference therein.}

We now characterize the {  SOS property and complete decomposability} for even order generalized Cauchy tensors with nonzero entries.
\bt\label{thm31} Let $\mathcal{C}$ be a generalized Cauchy tensor with even order $m$ and dimension $n$. Let $c, d\in \mathbb{R}^n$
be generating vectors of $\mathcal{C}$. Assume $d_i\neq0,~i\in [n]$. Then, the following statements are equivalent:

\begin{itemize}
\item[{\rm (i)}] the generalized Cauchy tensor $\mathcal{C}$ is a completely decomposable tensor;
\item[{\rm (ii)}] the generalized Cauchy tensor $\mathcal{C}$ is an SOS tensor;
\item[{\rm (iii)}] the generalized Cauchy tensor $\mathcal{C}$ is positive semi-definite;
\item[{\rm (iv)}] $c>{\bf 0}$.
\end{itemize}
\et
\proof Since $m$ is even, by the definitions of completely decomposable tensor, SOS tensor and positive
semi-definite tensor, we can easily obtain ${\rm (i)}\Rightarrow{\rm (ii)}$ and ${\rm (ii)}\Rightarrow{\rm (iii)}$.

$[{\rm (iii)} \Rightarrow {\rm (iv)}]$ ~ Let $\mathcal{C}$ be an even order generalized Cauchy tensor which is positive semi-definite.  Then
$$\mathcal{C} e_i^m=\frac{d_i^m}{m c_i} \ge 0.$$
So $c_i>0$ for all $i=1,\ldots,n$.

$[{\rm (iv)} \Rightarrow {\rm (i)}]$~Suppose that $c>{\bf 0}$. Then, for any $x\in \mathbb{R}^n$,
\begin{equation}\label{e32}
\begin{aligned}
f(x)=& \, \mathcal{C} x^m  = \sum_{i_1,i_2,\cdots,i_m=1}^n \frac{d_{i_1} d_{i_2} \cdots d_{i_m}}{c_{i_1}+c_{i_2}+\cdots+c_{i_m}} x_{i_1}x_{i_2} \cdots x_{i_m}  \\
= & \sum_{i_1,i_2,\cdots,i_m=1}^n  \left(\int_0^1  t^{c_{i_1}+c_{i_2}+\cdots+c_{i_m}-1} d_{i_1} d_{i_2} \cdots d_{i_m} x_{i_1}x_{i_2} \cdots x_{i_m} dt \right) \\
=&\int_0^1  \left(\sum_{i_1,i_2,\cdots,i_m=1}^n t^{c_{i_1}+c_{i_2}+\cdots+c_{i_m}-1} d_{i_1} d_{i_2} \cdots d_{i_m} x_{i_1}x_{i_2} \cdots x_{i_m} \right)dt  \\
=& \int_0^1  \left(\sum_{i=1}^n t^{c_i-\frac{1}{m}} d_ix_i
\right)^m dt.
\end{aligned}
\end{equation}
By the definition of Riemann integral, we have
$$
\mathcal{C} x^m = \lim_{k \rightarrow \infty}\sum_{j=1}^k
\frac{\left(\sum_{i=1}^n (\frac{j}{k})^{c_i-\frac{1}{m}} d_ix_i
\right)^m}{k}.
$$
Let $\mathcal{C}_k$ be the symmetric tensor such that
\begin{equation}\label{e33}
\begin{aligned}
\mathcal{C}_k x^m=&\sum_{j=1}^k \frac{\left(\sum_{i=1}^n (\frac{j}{k})^{c_i-\frac{1}{m}} d_ix_i \right)^m}{k} \\
= & \sum_{j=1}^k \left(\sum_{i=1}^n \frac{(\frac{j}{k})^{c_i-\frac{1}{m}} d_i}{k^{\frac{1}{m}}}x_i \right)^m \\
=& \sum_{j=1}^k \left(\langle u^j, x \rangle \right)^m,
\end{aligned}
\end{equation}
where
\begin{equation}\label{e34}
u^j=\left(\frac{(\frac{j}{k})^{c_1-\frac{1}{m}}
d_1}{k^{\frac{1}{m}}},\cdots, \frac{(\frac{j}{k})^{c_n-\frac{1}{m}}
d_n}{k^{\frac{1}{m}}}\right) \in \mathbb{R}^n, \ j=1,\cdots,k.
\end{equation}
{  Let ${\rm CD}_{m,n}$ denote the set consisting of all completely decomposable tensor with order $m$ and dimension $n$. From  \cite[Theorem 1]{LQX14}, ${\rm CD}_{m,n}$ is a closed convex cone when $m$ is even.} It then follows that
$\mathcal{C}=\lim_{k \rightarrow \infty}\mathcal{C}_k$ is also a  completely decomposable tensor.
\qed

{  Next, we provide a sufficient and necessary condition for the complete positivity of a generalized Cauchy tensor with nonzero entries, in terms of its generating vectors.}

\bt\label{thm32} Let $\mathcal{C}$  be a generalized Cauchy tensor defined as in Theorem \ref{thm31} with generating vectors $c,~d\in \mathbb{R}^n$.
Assume $d_i\neq0,~i\in [n]$.
Then $\mathcal{C}$ is a completely positive tensor
if and only if $c>0$ and $d >0$.
\et

\proof
For necessary condition, suppose that $\mathcal{C}$ is a completely positive tensor. Then, for any vector $x \in \mathbb{R}^n_+$, we must have
$\mathcal{C}x^m \ge 0$. So,
$\mathcal{C} e_i^m=\frac{d_i^m}{m c_i} \ge 0$.
This implies that $c>0$.
To finish the proof, we only need to show $d>0$. To see this, we proceed by the method of contradiction and suppose that
$$
I_{-}:=\left\{i \in \{1,\cdots,n\}: d_i<0\right\} \neq \emptyset.
$$
Denote $r$ to be the cardinality of $I_{-}$. Without loss of generality, we assume that $I_{-}=\{1,\cdots,r\}$.  Then,  $d_1<0$ and $d_{r+1}>0$, and hence, the $(r+1,1,\ldots,1)^{\rm th}$ entry of $\mathcal{C}$ satisfies
$$
\mathcal{C}_{r+1\cdots 11} = \frac{d_{r+1} d_1^{m-1}}{c_{r+1}+(m-1)c_1}<0.
$$
Note that each entry of a completely positive tensor   must be a nonnegative number.
This makes contradiction, and hence, the necessary condition follows.

To prove the sufficient condition, from (\ref{e32})-(\ref{e33}), we know that
$$
\mathcal{C}x^m=\lim_{k\rightarrow\infty}\sum_{j=1}^k \left(\langle u^j, x \rangle \right)^m.
$$
By conditions $c>0$, $d >0$ and by (\ref{e34}), we see that $u^j \in \mathbb{R}^n_+,~j\in[k]$.
So each $\mathcal{C}_k$ is a completely positive tensor.
{  Let ${\rm CP}_{m,n}$ denote the set consisting of all completely positive tensors with order $m$ and dimension $n$.
From \cite{QXX}, ${\rm CP}_{m,n}$ is a closed convex cone for any $m,n \in \mathbb{N}$.  It then follows that
$\mathcal{C}=\lim_{k \rightarrow \infty}\mathcal{C}_k$ is also a completely positive tensor.}
\qed

As a direct corollary of Theorem \ref{thm31} and \ref{thm32}, we obtain the following characterization of sum-of-squares property, positive semi-definiteness and complete positivity of Cauchy tensors. We note that the equivalence of {\rm (ii)} and {\rm (iv)} of this corollary was established in \cite{chen14}.
\begin{corollary}\label{corol31}
Let $c \in \mathbb{R}^n$. Let $\mathcal{C}$ be a Cauchy tensor as defined in Definition \ref{def31} with even order $m$.
 Then, the following statements are equivalent.
 \begin{itemize}
\item[{\rm (i)}] the  Cauchy tensor $\mathcal{C}$ is a SOS tensor;
\item[{\rm (ii)}] the  Cauchy tensor $\mathcal{C}$ is positive semi-definite;
\item[{\rm (iii)}] the  Cauchy tensor $\mathcal{C}$ is a completely positive tensor;
\item[{\rm (iv)}] $c>0$.
\end{itemize}
\end{corollary}

Let $\mathcal{A}=(a_{i_1\cdots i_m})$ and $\mathcal{B}=(b_{i_1\cdots i_m})$ be two real tensors with order $m$ and dimension $n$.
Then their Hadamard product is a real order $m$ dimension $n$
tensor
$$\mathcal{A}\circ \mathcal{B}=(a_{i_1\cdots i_m}b_{i_1\cdots i_m}).$$
From Proposition 1 of \cite{QXX}, we know that the Hadamard product of two completely positive tensors is also a completely positive tensor.
So, by Corollary \ref{corol31}, we have the following conclusion.
\begin{corollary}\label{corol32} Let $\mathcal{C}_1$ and $\mathcal{C}_2$ be two positive semi-definite Cauchy tensors. Then the Hadamard product
$\mathcal{C}_1 \circ \mathcal{C}_2$ is also positive semi-definite.
\end{corollary}

Next, we have the following theorem on H-eigenvalues of nonnegative Cauchy tensors. By \cite{Yang10}, we know that {  each nonnegative tensor
has at least one H-eigenvalue}, which is the largest modulus of its eigenvalues. Here, for nonnegative Cauchy tensors, all the $H$-eigenvalues
must be nonnegative.

\bt\label{thm33} Let  $\mathcal{C}$ be a nonnegative Cauchy tensor with order $m$ dimension $n$. Let
$c=(c_1,c_2,\cdots,c_n)^T$ be the generating vector of tensor $\mathcal{C}$. Then all H-eigenvalues of Cauchy tensor $\mathcal{C}$
are nonnegative.
\et
\proof In the case where $m$ is even, since $\mathcal{C}$ is nonnegative and the definition of a Cauchy tensor, we have $c_i>0$, $i=1,\ldots,n$.
From Theorem 2.1 of \cite{chen14}, we know that $\mathcal{C}$ is positive semi-definite. Then, Theorem 5 of \cite{Qi05} gives us
that all H-eigenvalues of $\mathcal{C}$ are nonnegative.

We now consider the case where $m$ is odd.  {  Let $\lambda$ be an} arbitrary H-eigenvalue of $\mathcal{C}$ with an H-eigenvector $x\neq0$.
By the definition of $H$-eigenvalue, it holds that
$$
\begin{aligned}
\lambda x_i^{m-1}=&(\mathcal{C} x^{m-1})_i \\
=& \sum_{i_2,\cdots,i_m=1}^n \frac{x_{i_2}x_{i_3} \cdots x_{i_m}}{c_{i}+c_{i_2}+\cdots+c_{i_m}}   \\
= & \sum_{i_2,\cdots,i_m=1}^n  \left(\int_0^1  t^{c_i+c_{i_2}+\cdots+c_{i_m}-1} x_{i_2}x_{i_3} \cdots x_{i_m} dt \right) \\
=&\int_0^1  \left(\sum_{i_2,\cdots,i_m=1}^n t^{c_{i}+c_{i_2}+\cdots+c_{i_m}-1} x_{i_2}x_{i_3} \cdots x_{i_m} \right)dt  \\
=& \int_0^1  \left(\sum_{j=1}^n t^{c_j+\frac{c_i-1}{m-1}} x_j\right)^{m-1} dt.
\end{aligned}
$$
This implies that $\lambda\geq0$ since $m$ is odd. Thus, the desired result holds. \qed

\setcounter{equation}{0}
\section{Further properties on Hankel tensors}
Hankel tensors arise from signal processing and some other applications \cite{BB, Ding14, Papy05, Qi14}. Recall that an order $m$ dimension $n$ tensor
$\mathcal{A}=(a_{i_1i_2\cdots i_m})$ is called a Hankel tensor if there is a vector
$v=(v_0, v_1,\cdots, v_{(n-1)m})^T$ such that
\begin{equation}\label{e41}
a_{i_1i_2\cdots i_m}=v_{i_1+i_2+\cdots +i_m-m},~\forall ~i_1,i_2,\cdots,i_m\in [n].\end{equation}
Such a vector $v$ is called the generating vector of Hankel tensor $\mathcal{A}$.

For any $k\in \mathbb{N}$, let $s(k,m,n)$ be the number of distinct sets of indices
$(i_1,i_2,\cdots,i_m)$, $i_j\in [n], j\in [m]$ such that $i_1+i_2+\cdots +i_m-m=k$. For example,
$s(0,m,n)=1, s(1,m,n)=m, s(2,m,n)=\frac{m(m+1)}{2}$. Suppose $\mathcal{P}=(p_{i_1i_2\cdots i_{(n-1)m}})$ is an order $(n-1)m$
dimension 2 tensor defined by
$$p_{i_1i_2\cdots i_{(n-1)m}}=\frac{s(k,m,n)v_k}{\binom {(n-1)m} k},$$
where $k=i_1+i_2+\cdots +i_{(n-1)m}-(n-1)m$. Then tensor $\mathcal{P}$ is called the
{\bf associated plane tensor} of Hankel tensor $\mathcal{A}$. When $n=2$, it is obvious that $\mathcal{P}=\mathcal{A}$.

In \cite{Qi14}, it was proved that, if a Hankel tensor is copositive, then its associated plane tensor $\mathcal{P}$
is copositive and the associated plane tensor is positive semi-definite if the Hankel tensor is
positive semi-definite. Since the associated plane tensor $\mathcal{P}$ has dimension 2, we can use the $Z$-eigenvalue method in \cite{QWW}
to check its positive semi-definiteness {  (alternatively, noting that any $2$-dimensional symmetric tensor is positive semi-definite if and only if it is a sums-of-squares tensor, we
can also verify the positive semi-definiteness of the associated plane tensor by {  solving} a semi-definite programming problem). Thus, the positive semi-definiteness of the associated plane tensor is a checkable necessary conditions for the positive semi-definiteness of even order Hankel tensors
(see more discussion in \cite{Qi14}).}
{  This naturally raises the following questions: Can these necessary conditions  be also sufficient? If not, are there any concrete counter-examples?}

{  We first present a result stating that the positive semi-definiteness of the associated plane tensor is equivalent to the
Vandermonde positive semi-definiteness of the original Hankel tensor.}

\bt\label{thm41}
{  Let $\mathcal{A}$ be a Hankel tensor defined as in (\ref{e41}) with an even order $m$. Then, the associated plane tensor
$\mathcal{P}$ is positive semi-definite if and only if $\mathcal{A}$ is Vandermonde positive semi-definite.}
\et
\proof For necessary condition, let $u=(1,\mu,\mu^2,\cdots,\mu^{n-1})^T\in \mathbb{R}^n$ be an arbitrary
Vandermonde vector.
If $\mu=0$, then we have
\begin{equation}\label{e42}
\mathcal{A}u^m=\sum_{i_1,i_2,\cdots,i_m\in [n]}a_{i_1i_2\cdots i_m}u_{i_1}u_{i_2}\cdots u_{i_m}=v_0.
\end{equation}
{  By our assumption}, for $y=(1,0)^T\in \mathbb{R}^2$, it follows that
$$\mathcal{P}y^{(n-1)m}=\sum_{i_1,i_2,\cdots,i_{(n-1)m}\in [2]}p_{i_1i_2\cdots i_{(n-1)m}}y_{i_1}y_{i_2}\cdots y_{i_{(n-1)m}}=v_0\geq0.$$
Combining this with (\ref {e42}), we obtain
\begin{equation}\label{e43}
\mathcal{A}u^m\geq0.
\end{equation}
If $\mu\neq0$, there exist $y_1, y_2\in \mathbb{R} \backslash \{0\}$ such that $\mu=\frac{y_2}{y_1}$. Let $y=(y_1,y_2)^T\in \mathbb{R}^2$. Then, we have
\begin{equation}\nonumber
\begin{aligned}
\mathcal{P}y^{(n-1)m}=&\sum\limits_{i_1,i_2,\cdots,i_{(n-1)m}\in [2]}p_{i_1i_2\cdots i_{(n-1)m}}y_{i_1}y_{i_2}\cdots y_{i_{(n-1)m}}\\
=&y_1^{(n-1)m}\sum\limits_{k=0}^{(n-1)m} \binom {(n-1)m} k \frac{s(k,m,n)v_k}{\binom {(n-1)m} k}\mu^k\\
=&y_1^{(n-1)m}\mathcal{A}u^m\\
\geq&0.
\end{aligned}
\end{equation}
By (\ref{e43}) and the fact that $m$ is even, for all Vandermonde vectors $u\in \mathbb{R}^n$, it follows that
$$\mathcal{A}u^m\geq0,$$
which implies Hankel tensor $\mathcal{A}$ is Vandermonde positive semi-definite.

For sufficiency, let $y=(y_1,y_2)^T\in \mathbb{R}^2$. We now verify that $\mathcal{P}y^{(n-1)m} \ge 0$. To see this, we first
consider the case where $y_1\neq0$. In this case, let $u=(1,\mu,\mu^2,\cdots,\mu^{n-1})^T\in \mathbb{R}^n$,
where $\mu=\frac{y_2}{y_1}$. From the analysis above, we have
\begin{equation}\label{e44}
\mathcal{P}y^{(n-1)m}=y_1^{(n-1)m}\mathcal{A}u^m\geq0
\end{equation}
since $m$ is even and $\mathcal{A}$ is Vandermonde positive semi-definite.

{  On the other hand, if $y=(y_1,y_2)^T\in \mathbb{R}^2$ with $y_1=0$}, then we let $y_\epsilon=(\epsilon, y_2)^T\in \mathbb{R}^2$ and $u=(1,\mu,\mu^2,\cdots,\mu^{n-1})^T\in \mathbb{R}^n$,
where $\mu=\frac{y_2}{\epsilon}$, $\epsilon>0$.
By (\ref{e44}), we have
$$\mathcal{P}y_\epsilon^{(n-1)m}=\epsilon^{(n-1)m}\mathcal{A}u^m\geq0.$$
Combining this with the fact that {  $\epsilon \mapsto \mathcal{P}y_\epsilon^{(n-1)m}$ is a continuous mapping}, it follows that
$$\mathcal{P}y^{(n-1)m}=\lim_{\epsilon\rightarrow0}\mathcal{P}y_\epsilon^{(n-1)m}\geq0.$$
This then implies that plane tensor $\mathcal{P}$ is positive semi-definite and the desired result holds.
\qed


Below, we provide an example illustrating that a Hankel tensor which is Vandermonde positive semi-definite need not to be
positive semi-definite. This example together with Theorem \ref{thm41}, also implies that the positive
semi-definiteness of the associate plane tensor is not sufficient for positive semi-definiteness of the Hankel tensor.

\begin{example}\label{exam41}   Let $\mathcal{A}$ be a Hankel tensor with order $m=4$ and dimension $n=3$.
{  Let the generating vector of $\mathcal{A}$ be $v_0=1, v_1=-1, v_2=1$ and $v_3=v_4=\cdots=v_8=0$.} So, for any $u=(1,\mu,\mu^2)^T\in \mathbb{R}^3$,
$$
\begin{aligned}
\mathcal{A}u^4=&\sum\limits_{i_1,i_2,i_3,i_4\in [3]}a_{i_1i_2i_3i_4}u_{i_1}u_{i_2}u_{i_3}u_{i_4}\\
=&\sum\limits_{k=0}^{k=8}s(k,4,3)v_k\mu^k\\
=&v_0+4v_1\mu+10v_2\mu^2\\
=&1-4\mu+10\mu^2\geq0
\end{aligned}
$$
for all $\mu\in \mathbb{R}$. By Theorem \ref{thm41}, we know that the associated plane tensor $\mathcal{P}$
is positive semi-definite. {  We now verify that $\mathcal{A}$ is not positive semi-definite. To see this,}  let $x=(1,1,-1)^T$, then,
$$
\begin{aligned} \mathcal{A}x^4 
=&\sum\limits_{i_1,i_2,i_3,i_4\in [3]}v_{i_1+i_2+i_3+i_4-4}x_{i_1}x_{i_2}x_{i_3}x_{i_4}\\
=&v_0x_1^4+4v_1x_1^3x_2+v_2(6x_2^2x_1^2+4x_1^3x_3)\\
=&x_1^4-4x_1^3x_2+(6x_2^2x_1^2+4x_1^3x_3)\\
=&1-4+6-4=-1<0,
\end{aligned} $$
which implies that Hankel tensor $\mathcal{A}$ is not positive semi-definite. \qed
\end{example}

{  The following example shows that the the copositivity of the associated plane tensor is also not sufficient for the copositivity of the Hankel tensor, in general.}

\begin{example}\label{exam42} Let $\mathcal{A}$ be a Hankel tensor with order $m=4$ and dimension $n=3$. {  Let the generating vector of $\mathcal{A}$ be $v_0=1, v_1=-1, v_2=\frac{1}{2}$, $v_3=v_4=\cdots=v_8=0$.} Let $x=(1,\frac{1}{2},0)^T$.
Then, we have
$$\mathcal{A}x^4=-\frac{1}{4}<0,$$
which implies that Hankel tensor $\mathcal{A}$ is not copositive.
On the other hand, it holds that
$$\mathcal{A}u^4=1-4\mu+5\mu^2\geq0$$
for any Vandermonde vector $u=(1,\mu,\mu^2)^T\in \mathbb{R}^3.$
By Theorem \ref{thm41}, the associated plane tensor $\mathcal{P}$ is positive semi-definite.
Thus, $\mathcal{P}$ is copositive. \qed
\end{example}

A special class of Hankel tensor is the complete Hankel tensors. {  To recall the definition of a complete Hankel tensor, we note that, for a Hankel tensor $\mathcal{A}$ with order $m$ dimension $n$}, if
\begin{equation}\label{e45}
\mathcal{A}=\sum\limits_{k=1}^{r}\alpha_k(u_k)^m,
\end{equation}
where $\alpha_k\in \mathbb{R}$, $\alpha_k\neq0$, $u_k=(1,\mu_k,\mu_k^2,\cdots,\mu_k^{n-1})^T\in \mathbb{R}^n$,
$k=1,2,\cdots,r$, for some $\mu_i\neq\mu_j$ for $i\neq j$, {  then, we say $\mathcal{A}$ has a {\bf Vandermonde decomposition}. The corresponding vector $u_k$, $k=1,\ldots,r$ are called {\bf Vandermonde vectors} and the minimum value
of $r$ is called {\bf Vandermonde rank} of $\mathcal{A}$ \cite{Qi14}.} From Theorem 3 of \cite{Qi14}, we know that $\mathcal{A}$ is a Hankel tensor if and only if
it has a Vandermonde decomposition (\ref{e45}).
If $\alpha_k>0$ for $k\in [r]$ in (\ref{e45}), then $\mathcal{A}$ is called a complete Hankel tensor.

In \cite{Qi14}, it is proved that an even order complete Hankel tensor is positive semi-definite. Moreover, examples {were also presented to show that the converse is, in general, not true}.
{  Here, in the following theorem, we show that if the Vandermonde rank of a Hankel tensor $\mathcal{A}$ is less than the dimension of the underlying space, then positive semi-definiteness
of $\mathcal{A}$ is equivalent to the fact that $\mathcal{A}$ is a complete Hankel tensor, and so, is further equivalent to the SOS property of $\mathcal{A}$.}

{  \bt\label{thm42}  Let $\mathcal{A}$ be a Hankel tensor with an even order $m$. Assume that the
Hankel tensor $\mathcal{A}$ has Vandermonde decomposition (\ref{e45}) with the Vandermonde rank $r$ satisfies $r\leq n$. Then, the following statements
are equivalent:
\begin{itemize}
 \item[{\rm (i)}] $\mathcal{A}$ is a positive semi-definite tensor;
 \item[{\rm (ii)}] $\mathcal{A}$ is a complete Hankel tensor.
  \item[{\rm (iii)}] $\mathcal{A}$ is an SOS tensor;
\end{itemize}
\et

\proof  We first note that the implications $[{\rm (ii)}] \Rightarrow [{\rm (iii)}]$ and $[{\rm (iii)}] \Rightarrow [{\rm (i)}]$ are direct consequences from the
definitions. Thus, to see the conclusion, we only need to prove $[{\rm (i)}] \Rightarrow [{\rm (ii)}]$. To do this, we proceed by the method of contradiction and
 assume that there exists at least one coefficient $\alpha_i$ in (\ref{e45}) which is negative. }
 Without loss of generality, we assume that $\alpha_1<0$.
For any $x=(x_1,x_2,\cdots,x_n)^T\in \mathbb{R}^n$, then we have
\begin{equation}\label{e46}
\begin{aligned}
\mathcal{A}x^m=&\sum\limits_{k=1}^r\alpha_k(u_k^Tx)^m\\
=&\alpha_1(u_1^Tx)^m+\alpha_2(u_2^Tx)^m+\cdots+\alpha_r(u_r^Tx)^m.
\end{aligned}
\end{equation}
Consider the following two homogeneous linear equation systems
$$Ax={\bf 0},~~Bx={\bf 0},$$
where
$$
A=\left(
\begin{matrix}
1 & \mu_1 & \mu_1^2 & \cdots & \mu_1^{n-1} \\
1 & \mu_2 & \mu_2^2 & \cdots & \mu_2^{n-1} \\
\vdots & \vdots & \vdots & \vdots & \vdots \\
1 & \mu_r & \mu_r^2 & \cdots & \mu_r^{n-1} \\
\end{matrix}
\right),
\quad B=\left(
\begin{matrix}
1 & \mu_2 & \mu_2^2 & \cdots & \mu_2^{n-1} \\
1 & \mu_3 & \mu_3^2 & \cdots & \mu_3^{n-1} \\
\vdots & \vdots & \vdots & \vdots & \vdots \\
1 & \mu_r & \mu_r^2 & \cdots &\mu_r^{n-1} \\
\end{matrix}
\right).
$$
By conditions $r\leq n$, it is easy to get
$$Rank(A)=r\leq n,~~Rank(B)=r-1<n,$$
which imply that there is vector $x_0\in \mathbb{R}^n$, $x_0\neq {\bf 0}$ such that
$$Ax_0\neq {\bf 0},~~Bx_0= {\bf 0}.$$
Here, $Rank(A)$ denotes the rank of matrix $A$. So, it holds that
$$u_1^Tx_0\neq0,~~u_i^Tx_0=0,~~i\in \{2,3,\cdots,r\}.$$
Combining this with (\ref{e46}), we have
$$\mathcal{A}x_0^m=\alpha_1(u_1^Tx_0)^m<0$$
since $m$ is even. However, {  this contradicts to} the fact that $\mathcal{A}$ is positive semi-definite.
Thus, all coefficients in (\ref{e45}) are positive and $\mathcal{A}$ is a complete Hankel tensor.
\qed

{  An interesting consequence of Theorem \ref{thm42} is as follows: a necessary condition for a PNS (positive semi-definite but not sum-of-squares)  Hankel tensor $\mathcal{A}$
is that the Vandermonde rank $r$ of the Hankel tensor $\mathcal{A}$ is strictly larger than the dimension $n$ of the underlying space. We note that searching for a
PNS Hankel tensor is a non-trivial task and is related to Hilbert's 17th question. Recently, some extensive study for PNS Hankel tensor has been initialed in
\cite{LQW14}.

Next, we provide some necessary conditions for the positive semi-definiteness of a Hankel tensor $\mathcal{A}$ in terms of the sign properties of the coefficients of its Vandermonde decomposition.}

\bp\label{prop41} Let $\mathcal{A}$ be a Hankel tensor with the Vandermonde decomposition (\ref{e45}). Suppose that $\mathcal{A}$
is positive semi-definite. Then,

(i) the coefficients of the Vandermonde decomposition satisfy
$$\alpha_1+\alpha_2+\cdots+\alpha_r\geq0;$$

(ii) if $r>n$, then the total number of positive coefficients of the Vandermonde decomposition is greater than or equal to $n$;

(iii) if $r\leq n$, then all coefficients of the Vandermonde decomposition are positive.
\ep

\proof
(i) Since $\mathcal{A}$ is positive semi-definite, so we have
$$\mathcal{A}e_1^m=\sum\limits_{i=1}^r \alpha_i (u_i^Te_1)^m=\alpha_1+\alpha_2+\cdots+\alpha_r\geq0.$$

(ii) {  Denote the total number of positive coefficients in (\ref{e45}) by $t$}. Without loss of generality, let
$$\alpha_i>0,~i\in [t];~\alpha_j<0,~j\in \{t+1,t+2,\cdots, r\}.$$
{  We proceed by the method of contradiction and suppose that $t<n$.
If $t=0$, we can easily get a contradiction because $\mathcal{A}$ is positive semi-definite.} If $0<t<n$, consider the following two linear equation systems
\begin{equation} \label{e47} Ax={\bf 0}\end{equation}
and
\begin{equation} \label{e48} Bx={\bf 0},\end{equation}
where
$$
A=\left(
\begin{matrix}
1 & \mu_1 & \mu_1^2 & \cdots & \mu_1^{n-1} \\
1 & \mu_2 & \mu_2^2 & \cdots & \mu_2^{n-1} \\
\vdots & \vdots & \vdots & \vdots & \vdots \\
1 & \mu_t & \mu_t^2 & \cdots & \mu_t^{n-1} \\
\end{matrix}
\right),
\quad B=\left(
\begin{matrix}
1 & \mu_1 & \mu_1^2 & \cdots & \mu_1^{n-1} \\
1 & \mu_2 & \mu_2^2 & \cdots & \mu_2^{n-1} \\
\vdots & \vdots & \vdots & \vdots & \vdots \\
1 & \mu_r & \mu_r^2 & \cdots & \mu_r^{n-1} \\
\end{matrix}
\right).
$$
{  Noting that
$Rank(A)=t<n$ and $Rank(B)=n$, basic linear algebra theory implies that the system (\ref{e47}) has nonzero solutions and system (\ref{e48}) has only zero solution.}
Thus, there exists $\bar{x}\in \mathbb{R}^n, \bar{x}\neq 0$ such that
$$u_i^T\bar{x}=0,~~i\in [t]  {  \mbox{ and } (u_{t+1}^T\bar{x},\ldots,u_{r}^T \bar x)^T \neq {\bf 0}}.$$
Note that the order $m$ is an even number (as $\mathcal{A}$ is positive semi-definite). This implies that
$${  \mathcal{A}\bar{x}^m}=\sum\limits_{j=t+1}^r\alpha_j(u_j^T\bar{x})^m<0.$$
This contradicts with the fact that $\mathcal{A}$ is positive semi-definite. Then we get $t\geq n$.

(iii) If $r\leq n$, then the conclusion is a direct result of Theorem \ref{thm42}. \qed

%
%

\setcounter{equation}{0}
\section{Properties of Cauchy-Hankel tensors}

{  In the literature, there is an important class of structured matrices called Cauchy-Hankel matrices which is closely related with Cauchy matrices and Hankel matrices simultaneously \cite{GAD04, solak02, solak03}. A matrix $A$ is called a Cauchy-Hankel matrices if it can be
formulated as
$$A=\left( \frac{1}{g+h(i+j)}\right)_{i,j=1}^n,$$
where $g$ and $h$ are real constants such that $h\neq0$ and $\frac{g}{h}$ is not an integer \cite{Boz98}.}

As a natural extension of Cauchy-Hankel matrix, a tensor $\mathcal{A}=(a_{i_1i_2\cdots i_m})$ with order $m$ and dimension $n$ is called a {\bf Cauchy-Hankel tensor}
if
\begin{equation}\label{e51} a_{i_1i_2\cdots i_m}=\frac{1}{g+h(i_1+i_2+\cdots+i_m)},~i_j\in [n],~j\in [m],\end{equation}
where $g, h\neq0\in \mathbb{R}$ and $\frac{g}{h}$ is not an integer.

It is obvious that a Cauchy-Hankel tensor is a symmetric tensor. From Definition \ref{def31},
we know that a Cauchy-Hankel tensor defined by (\ref{e51}) is a Cauchy tensor \cite{chen14} with generating vector
$$c=(\frac{g}{m}+h,\frac{g}{m}+2h,\cdots,\frac{g}{m}+nh)^T\in \mathbb{R}^n,$$
and it is a Hankel tensor \cite{Papy05, Qi14} at the same time with
$$v_k=\frac{1}{g+h(k+m)},~k\in\{0,1,2,\cdots,(n-1)m\}.$$

\bt\label{thm51} {  Let $\mathcal{A}$ be a Cauchy-Hankel tensor defined as in (\ref{e51}) with an even order $m$.}
 Then, $\mathcal{A}$ is positive definite if and only if
$$g+mh>0,~~g+nmh>0.$$
\et

\proof For necessary condition, since $\mathcal{A}$ is positive definite, so we have
$$\mathcal{A}e_1^m=\frac{1}{g+mh}>0,~~\mathcal{A}e_n^m=\frac{1}{g+mnh}>0,$$
and the desired results hold.

For sufficiency, since
$$g+mh>0,~~g+nmh>0,$$
it follows that
$$g+smh>0,~~\forall~s\in \{1,2,\cdots,n\}.$$
Combining Theorem 2.3 of \cite{chen14} and the fact that
$$g+imh\neq g+jmh,~~\forall~i,j\in [n],~{  i\neq j},$$
 we know that $\mathcal{A}$ is positive definite and
the desired result follows.
\qed

Next, we define the homogeneous polynomial $f(x)$ as below
$$f(x)=\mathcal{A}x^m =\sum\limits_{i_1,i_2,\cdots,i_m\in [n]}a_{i_1i_2\cdots i_m}x_{i_1}x_{i_2}\cdots x_{i_m},
$$
for $x=(x_1,x_2,\cdots,x_n)^T\in \mathbb{R}^n$. Let $x,y\in X \subseteq \mathbb{R}^n$. If $f(x)\geq f(y)$ for
any $x\geq y (x\leq y$ respectively), then  we say $f(x)$ is monotonically increasing (monotonically decreasing respectively) in $X$.
If $f(x)> f(y)$ for
any $x\geq y, x\neq y (x\leq y, x\neq y$ respectively), then  we say $f(x)$ is strict monotonically increasing (strict monotonically decreasing respectively) in $X$.

When $\mathcal{A}$ is a Cauchy tensor with even order, it has been proved that $f(x)$ is strict monotonically increasing in $\mathbb{R}^n_+$ if
the Cauchy tensor $\mathcal{A}$ is positive definite; { while the converse need not to be true \cite{chen14}}.  For even order Cauchy-Hankel tensors, we have the following conclusion, which is stronger than the
corresponded conclusion listed in \cite{chen14}.

\bt\label{thm52} {  Let $\mathcal{A}$ be a Cauchy-Hankel tensor defined as in (\ref{e51}) with an even order $m$}. Then, $\mathcal{A}$
is positive definite if and only if $f(x)=\mathcal{A}x^m$ is strict monotonically increasing in $\mathbb{R}^n_+$.
\et
\proof For the only if part, suppose $x,y\in \mathbb{R}^n_+$, $x\geq y$ and $x\neq y$, which means that there exists at least
one subscript $i$ satisfying $x_i>y_i$. Then, we have
$$
\begin{aligned} f(x)-f(y)=&\mathcal{A}x^m-\mathcal{A}y^m\\
=&\sum\limits_{i_1,i_2,\cdots,i_m\in [n]}\frac{x_{i_1}x_{i_2}\cdots x_{i_m}-y_{i_1}y_{i_2}\cdots y_{i_m}}{g+h(i_1+i_2+\cdots+i_m)} \\
=&\frac{x_i^m-y_i^m}{g+imh}+\sum\limits_{i_1i_2\cdots i_m\neq ii\cdots i}\frac{x_{i_1}x_{i_2}\cdots x_{i_m}-y_{i_1}y_{i_2}\cdots y_{i_m}}{g+h(i_1+i_2+\cdots+i_m)}.\\
\end{aligned}
$$
Since $\mathcal{A}$ is positive definite, by Theorem \ref{thm51}, we obtain
$$g+kmh>0,~~\forall ~k\in [n].$$
So, we obtain
$$\frac{x_i^m-y_i^m}{g+imh}>0$$
and
$$\sum\limits_{i_1i_2\cdots i_m\neq ii\cdots i}\frac{x_{i_1}x_{i_2}\cdots x_{i_m}-y_{i_1}y_{i_2}\cdots y_{i_m}}{g+h(i_1+i_2+\cdots+i_m)} \geq 0.$$
Thus, we have
$$f(x)-f(y)>0,$$
which implies that $f(x)$ is strict monotonically increasing in $\mathbb{R}^n_+$.

For the if part, note that $e_i\in \mathbb{R}^n_+$ and $e_i\geq {\bf 0}$, $e_i\neq {\bf 0}$, $i=1,n$.
{  It then follows that }
$$f(e_1)-f({\bf 0})=\mathcal{A}e_1^m=\frac{1}{g+mh}>0$$
and
$$f(e_n)-f({\bf 0})=\mathcal{A}e_n^m=\frac{1}{g+nmh}>0.$$
By Theorem \ref{thm51}, we know that Cauchy-Hankel tensor $\mathcal{A}$ is positive definite and the desired results hold.
\qed

\section{Final Remarks}
In this article, we present various new results on Cauchy tensors and Hankel tensors which complements the existing literature.
Firstly, we show that positive semi-definite Cauchy tensors are SOS tensors.
Furthermore, we prove that an even order Cauchy tensor is positive semi-definite if and only if it is a completely
positive tensor. The nonnegativity of H-eigenvalues of nonnegative Cauchy tensors are also established.
{  For Hankel tensors, we  prove that it is Vandermonde positive semi-definite if and only if the associated plane tensor is positive semi-definite.  We also show that,
if the Vandermonde rank of a Hankel tensor $\mathcal{A}$ is less than the dimension of the underlying space, then positive semi-definiteness
of $\mathcal{A}$ is equivalent to the fact that $\mathcal{A}$ is a complete Hankel tensor, and so, is further equivalent to the SOS property of $\mathcal{A}$.}
Finally, properties of Cauchy-Hankel are also given.

\end{document}